# Designing a Circular Economy Network for PPE Masks Supply Chain: A Case Study of British Columbia, Canada

Jainil Dharmil Shah[1,2[0000-0007-2499-204X]], Behrooz Khorshidvand[1[0000-0001-5597-5188]], Niloofar Gilani Larimi[1[0000-0002-3823-7639]], Adel Guitouni[1[0000-0003-2336-5398]]

[1] Peter B. Gustavson School of Business, University of Victoria, Victoria, BC, Canada
[2] Edwardson School of Industrial Engineering, Purdue University, West Lafayette, IN, USA

**Abstract.** In recent years, there has been growing interest in building closed-loop supply chains (SCs). However, many of the current methods struggle when it comes to fully embracing circular economy principles. Enhancing the design and management of these networks holds significant potential to promote stronger collaboration among supply chain partners, ultimately fostering more sustainable and efficient operational practices. For this purpose, this study addresses a new circular economy model based on a real case study of a mask SC in the healthcare sector in British Columbia, Canada. The objective is to show that implementing circular practices can lead to considerable financial and environmental benefits, as well as the creation of new job opportunities. To achieve this, a multi-objective mixed-integer linear programming model is developed to identify the most efficient trade-off among sustainable objectives, while adhering to imposed constraints. The proposed closed-loop SC model outperforms the existing linear model in all three aspects of sustainability, namely economic, environmental, and social. The improvement leads to significant economic and environmental benefits by preventing the disposal of used masks, giving them a second chance for disinfection and reprocessing, and reintroducing them into the SC as new ones. While the results of the circular economy model demonstrate profit gains and environmental recovery, the linear SC model showed negative profit with higher carbon emissions. Regarding the social aspects, compared to the current system, our approach not only nearly doubles the number of jobs created but also significantly reduces shortages, highlighting sustainable development aspects related to equity and social welfare. The findings offer valuable insights for researchers and practitioners seeking to implement sustainability within a circular economy framework in SCs.

## 1 Introduction

Personal Protective Equipment (PPE) includes essential gear such as N95 masks, surgical masks, gowns, and goggles designed to prevent infection or injury to the wearer. The majority of the inputs and raw materials used to make PPE are contracted out to low-cost suppliers. The production of PPE depends mainly on low-cost suppliers for raw materials such as polyester, polyamide, and cotton fiber. These materials are typically sourced from various countries around the world. Once acquired, manufacturers process them into protective clothing, which is then sold to consumers.

The COVID-19 pandemic caused an unprecedented surge in demand for PPE, overwhelming existing production capacities. Production and logistical limitations hindered the supply chain's ability to meet this increased demand, resulting in raw material shortages and a global backlog of supply orders spanning four to six months [1]. Many countries-imposed export bans on PPE and key raw materials, exacerbating the scarcity. Consequently, the cost of PPE increased dramatically: surgical masks surged sixfold, respirators trebled, and gowns doubled in price [2].

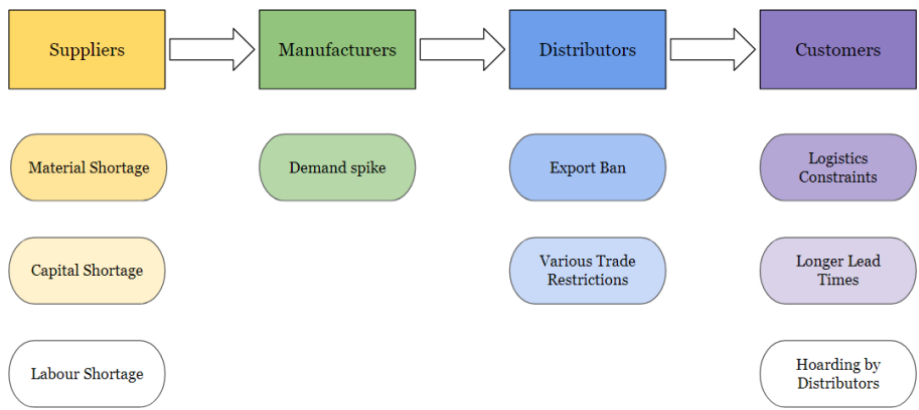

**Fig. 1.** PPE Supply Chain Bottlenecks

During significant public health crises like the COVID-19 pandemic, hospitals generate increased medical waste, posing severe environmental and health risks [3]. Regulations and policies have been implemented to manage medical waste effectively, but most techniques follow a conventional linear economy model characterized by a "take-make-waste" approach. This unsustainable model depletes natural resources and generates substantial waste and harmful emissions, contributing to environmental degradation and climate change [4].

During the pandemic, the sharp rise in demand for personal protective equipment (PPE) brought new attention to the issue of waste. Most PPE are designed for single use, which makes proper disposal a major challenge. A circular economy approach, centered on reusing, reprocessing, and recycling materials, offers a more sustainable alternative to this linear model of consumption and disposal [5- 6]. According to a study conducted by Vitacore [7], over 63,000 tons of single-use masks and respirators related to Covid-19 were projected to be used in Canada over the AY'22. In December 2020, estimates indicated that approximately 1.5 billion masks entered the ocean in 2020 alone [8]. According to a study by Evreka [9], a single face mask can release 173,000 microfibers into our oceans and seas.

In addition, if only 1% of the masks are disposed of incorrectly, it would result in 10 million masks per month. Calculating the numbers, each mask weighs approximately 4 grams, amounting to over 88,000 lbs. of plastic entering the environment each month [10]. These masks are typically made from plastics such as polypropylene and do not biodegrade. It is projected that this plastic will take around 400 years to disintegrate, which means that 15 generations later, the masks will be nearly gone if we stop adding them to the environment today [11]. Another significant issue that Canada faced during the Covid-19 pandemic which can also arise in any disastrous situation was the uncertain supply and demand of personal protective equipment (PPE), particularly in British Columbia (BC). According to the data on one month's PPE consumption by the different healthcare centres/hospitals under various regional health authorities in BC, the number of available international suppliers dropped to zero, necessitating healthcare decision-makers in BC to obtain the necessary PPEs from other regions, primarily Ontario. Given these challenges, this research aims to address two significant questions:

1. **Comparative Analysis**: Is a circular economy model superior to a linear economy model for the PPE supply chain with particular focus on masks in British Columbia, considering the three pillars of sustainability: economic, environmental, and social?

2. **Optimization**: How can we optimize the circular economy model regarding the location of collection, reprocessing, and recycling facilities to minimize costs, reduce environmental impact (carbon footprint), and maximize job creation?

It considers economic concerns while striving to improve environmental and social outcomes, incorporating uncertainties in both demand and supply. The study aims to design and optimize a multi-objective circular economy model for the mask supply chain, integrating economic, environmental, and social factors into the multi-objective function.

The structure of the paper is as follows: Section 2 provides a review of the existing literature relevant to our study. Section 3 defines the problem and presents the proposed circular economy framework. Section 4 details the objective functions and constraints of our mathematical model. Section 5 discusses the results, and Section 6 concludes with future research directions.

## 2     Literature Review

This study is closely related to four streams of literature on supply chain management: circular economy, sustainable supply chain, healthcare supply chain, and waste management.

The issue of medical waste and its environmental impact had become a pressing concern, especially in the context of the COVID-19 pandemic. As noted by [12], medical waste includes any material that has come in contact with patients, healthcare staff, or related workers during treatment or care. Due to its potentially infectious or hazardous nature, this type of waste poses serious risks- not only to human health but also to the environment- contributing to both micro- and macro-plastic pollution. The burden of environmental contaminants on public health is considerable. The healthcare sector is estimated to be responsible for around 4.6 percent of global greenhouse gas emissions. A 2020 estimate indicated that the global carbon footprint of healthcare accounts for 4.4% of all greenhouse gas emissions, despite health spending constituting 10% of global GDP [13].

Healthcare involves the continuous maintenance of health through the prevention, diagnosis, and treatment of diseases and addressing physical and mental disabilities. Timely action, precision, and positive outcomes are paramount. The management and coordination of logistics and supply chains for medical equipment and supplies are referred to as the healthcare supply chain [14]. Healthcare Supply Chain Logistics encompasses processes and workforce activities across different teams, facilitating the movement of medicines, surgical equipment and other necessary products. These logistics include pharmaceutical products, medical & surgical supplies, devices and other items required by healthcare professionals such as doctors, nurses, and administrative staff [15]. The complexity and unpredictability of supply chains present challenges for the healthcare sector [16]. Information and communication technology (ICT) and contract redesign are two methods to improve inefficiencies in healthcare supply chains [17]. Unlike conventional supply chains, where the main objective is maximizing profit, healthcare supply chains focus on optimizing efficiency and effectiveness of treatment. Efficient healthcare supply chains lead to improved processes, better resource utilization, satisfied employees, effective treatment and happy patients.

The transition to single-use medical equipment has harmed the environment and public health, strained supply chains and increased healthcare costs. Single-use disposables epitomize a linear economy (or "take-make-waste") where products are manufactured, used once, and disposed of [18]. Conversely, a circular economy is restorative or regenerative by design, aiming to reduce resource input, waste, emissions, and energy leakage by closing material and energy loops. This involves improving product longevity and recovering value from waste through industrial symbiosis, where one industry's outputs serve as raw materials for another. Industrial ecology, biomimicry, natural capitalism and the performance economy have all contributed ideas to the circular economy concept [19].

Sustainable supply chains integrate and coordinate the flow of goods and services within and among businesses, focusing on social responsibility and long-term profitability [20]. Sustainable supply chain management can result in cost savings, improved brand reputation and new business opportunities [21]. Organizations must overcome challenges such as drivers, facilitators, barriers and impact assessment factors to adopt sustainable supply chain management practices [22]. Managers must understand sustainability-related challenges and develop collaborative supply chain strategies [23]. Sustainable development in healthcare supply chain networks is essential to address resource utilization, waste management, and inclusive economic growth [24]. When healthcare systems adopt more sustainable approaches, they not only boost operational efficiency but also help move the industry closer to the principles of a circular economy [25]. However, a range of challenges, including financial, technological, organizational, policy-related, as well as social and cultural can make this transition difficult [26]. Overcoming these obstacles calls for smarter inventory management, more efficient processes and strategies that help reduce overall costs. Integrating supply chain practices, managing supply chain risks and adopting supply chain 4.0 can enhance healthcare supply chain performance and support the circular economy.

The Ellen MacArthur Foundation defines the circular economy as "an industrial economy that is restorative or regenerative by intention and design" [27]. The concept of the circular economy has evolved as societies and businesses recognize the need for sustainable resource management and waste reduction.

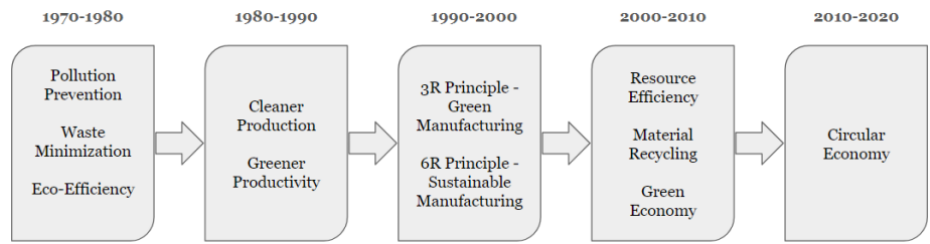

**Fig. 2.** Evolution of the concept of Circular Economy

The transition from the healthcare industry's linear paradigm to a circular economy model holds great promise, particularly in the context of the ongoing demand for PPE and its environmental impacts. The strategies proposed by the Canadian Coalition for Green Health Care [28], based on the three elements of the waste hierarchy—Reduce, Reuse, and Recycle—offer a hopeful path toward creating a closed-loop supply chain for different components of PPE.

In [29], a data-driven digital transformation approach was explored to optimize reverse logistics in medical waste management systems, particularly in the wake of the COVID-19 pandemic. This study highlighted that data-driven digital transformation can significantly enhance the efficiency of medical waste management processes. Decision makers in different echelons of a supply chain can optimize their operational decisions using real-time data analytics, which leads to reduced costs and improved service delivery. Regarding PPE waste management, [30] emphasized the importance of the 3R principles (i.e., reduce, reuse, and recycle). They highlighted that recycling

PPE benefits both environmental aspects and profitability. However, these studies currently focus only on economic & environmental objective functions and do not capture the social pillar of sustainability influencing reverse logistics operations.

Our study underscores the critical importance of integrating circular economy principles into the PPE supply chain. Doing so is not just a matter of environmental, economic, and social sustainability, but a fundamental necessity for the future of healthcare.

## 3 Problem Definition

British Columbia has five regional health authorities: Northern Health Authority, Island Health Authority, Vancouver Island Health Authority, Vancouver Coastal Health Authority, and Fraser Health Authority. These regional health authorities operate under the Provincial Health Services Authority (PHSA), which manages and regulates their activities.

Currently, the PPE supply chain in British Columbia follows a linear economy model comprising six echelons: suppliers, the PHSA, five regional health authorities (RHAs), various healthcare departments, and disposal centers. This traditional "take-make-dispose" model, depicted in Figure 3, represents the prevailing approach in which products are manufactured, used once, and discarded.

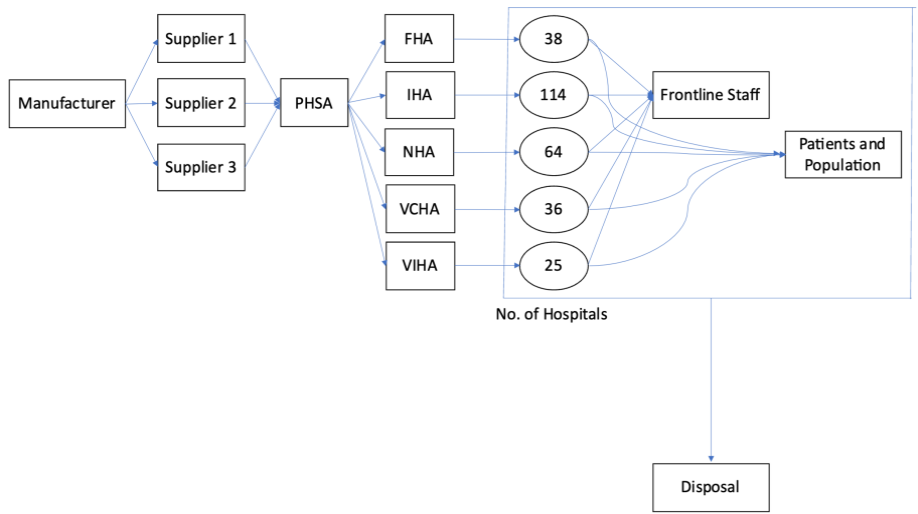

**Fig. 3.** Linear Economy Model for British Columbia

The COVID-19 pandemic highlighted significant vulnerabilities in this linear supply chain. The sudden and dramatic increase in demand for PPE overwhelmed existing production capacities, revealing critical limitations in production and logistics. As a result, the supply chain struggled to meet the heightened demand, leading to severe shortages and a backlog of supply orders that extended globally. Moreover, many countries implemented export bans on PPE and key raw materials, further complicating the situation. This scarcity led to a sharp rise in PPE costs, placing additional strain on healthcare systems striving to secure essential supplies.

For this study, Vancouver Island Health Authority (VIHA) was chosen as a region of focus to create a detailed model that can later be extrapolated to other health authorities. Within VIHA, the linear economy model includes 25 healthcare centers (hospitals), as illustrated in Figure 4. This focused approach allows for a thorough analysis of this regional authority's specific challenges and opportunities.

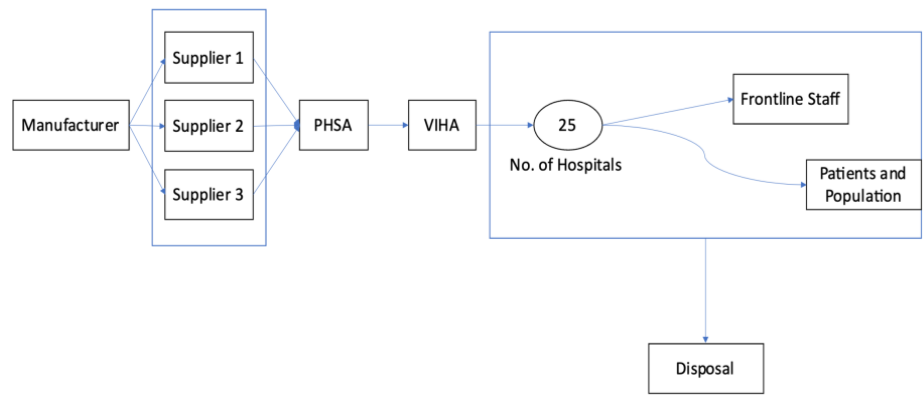

**Fig. 4.** Linear Economy Model for Vancouver Island Health Authority

The linear economy leads to significant resource depletion and waste generation. PPE items, once used, are discarded, contributing to mounting environmental pollution and medical waste accumulation. The linear model is highly vulnerable to disruptions. The dependency on global suppliers for raw materials and finished products makes the supply chain susceptible to geopolitical issues and logistical challenges, as seen during the pandemic. Moreover, the linear approach incurs substantial costs due to the procurement and transportation of raw materials and waste management and disposal processes.

A Closed Loop Supply Chain Network (CLSCN) based on the circular economy model is proposed to address these challenges. This model includes eight echelons: the original six from the linear economy model, additional collection centers, and reprocessing and recycling centers, as shown in Figure 5. Used masks collected from frontline staff, patients, and the common public are gathered at designated collection centers. At these centers, each mask is carefully inspected to assess its condition. Based on this assessment, a decision is made whether the mask should be sent for disinfection followed by reprocessing, or if it is unsuitable for reuse, sent directly for disposal. Masks that meet the criteria for reprocessing undergo thorough sterilization procedures before being repurposed, helping to reduce waste and extend the life cycle of PPE. This approach not only minimizes environmental impact but also supports the circular economy by reintegrating reusable materials into the supply chain.

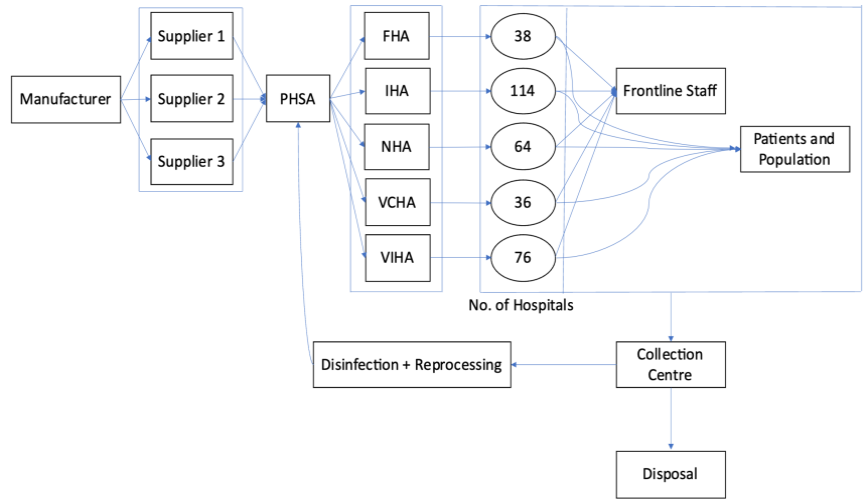

**Fig. 5.** Proposed Circular Economy Model for British Columbia

For our study, recycling centres for masks have been excluded and instead focus has been laid on reprocessing and disinfection facilities, which provide more excellent value for the healthcare sector. 25 potential locations in almost every majorly populated area of Vancouver Island Health Authority have been considered for setting up either a collection facility or a disinfection and reprocessing facility. The distance matrix between these 25 potential locations for setting up of collection and reprocessing facilities is calculated using GEOPY library in Python. Used masks by frontline staff, patients and general population in the surrounding regions of these hospitals/healthcare centres will be collected at collection centres where a decision would be taken on whether these masks are fit for reuse or need to be sent for directly disposal. The ones deemed fit for reuse would be sent to the reprocessing centres where they would be sterilized/disinfected using an appropriate decontamination technique and would be sent back to PHSA.

The designed and optimized circular economy model for the Vancouver Island Health Authority data is depicted in Figure 6.

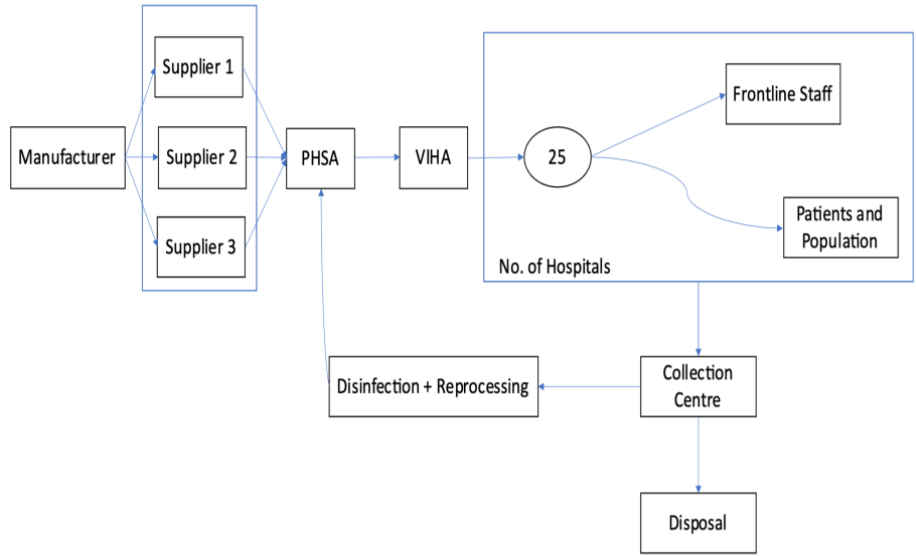

**Fig. 6.** Proposed Circular Economy Model for our study of Vancouver Island Health Authority

In this study, we consider a multi-objectives model for the Circular Economy Model:

1. **Economic objective:** The first objective is minimizing the total costs incurred throughout the mask supply chain. This includes optimizing collection locations, reprocessing, and disinfection facilities to reduce transportation and operational costs. The circular economy model can also reduce procurement costs by reusing materials and decreasing disposal costs.

2. **Environmental objective:** The second objective is to reduce the environmental impact of critical processes within the mask and respirator supply chain. By extending the lifecycle of PPE through reuse and reprocessing, the model reduces the need for new raw materials and lowers greenhouse gas emissions. Environmental impact data will be obtained from Life Cycle Assessment (LCA) results, measured in terms of $CO_2$ equivalent emissions.

3. **Social objective:** The third objective is to maximize job creation within the reverse logistics network of the mask supply chain. Establishing collection and reprocessing facilities generates employment opportunities, contributing to local economic development and social sustainability.

- **Economic objective:** Our primary economic goal is to minimize the overall costs associated with the PPE supply chain. This involves careful selection of locations for collection, reprocessing, and disinfection facilities to optimize transportation routes and reduce operational expenses. By integrating reuse and reprocessing into the supply chain, we aim to lower procurement costs and diminish the financial burden of disposal.

- **Environmental objective:** The environmental goal focuses on minimizing the carbon footprint and other adverse environmental impacts of PPE production and disposal. Implementing a circular economy model can significantly reduce greenhouse gas emissions by decreasing the demand for virgin raw materials. Life Cycle Assessment (LCA) data will be utilized to quantify these environmental benefits, ensuring that the model's impact is comprehensively evaluated.

- **Social objective:** Our social objective emphasizes enhancing community welfare through job creation. By establishing new facilities for collection and reprocessing, we aim to generate employment opportunities that contribute to local economic growth and social well-being. This approach supports the broader sustainable development goal by integrating economic, environmental, and social benefits.

By transitioning from a linear to a circular economy model, we aim to create a more resilient, sustainable, and cost-effective PPE supply chain for the Vancouver Island Health Authority. This approach can serve as a blueprint for other regions and healthcare systems seeking to enhance their supply chain sustainability.

## 3.1 Assumptions and Limitations of the Proposed Model

1. Focus is on only one critical component of PPE, i.e. masks, in particular all types of FFP2 respirators (N95 masks, etc) which was highest in demand and extensively used during the entire Covid-19 pandemic.

2. Usage of all types of FFP2 respirators at hospitals has been modelled and simulated randomly based on one month's real-time data on PPE consumption in British Columbia from September 1, 2020, to September 30, 2020.

3. For ease of network design and modelling, usage of PPE has been considered at and in areas which are in the vicinity of 25 major hospitals/healthcare centre locations within Vancouver Island Health Authority (VIHA), possibility of existence of other hospitals/healthcare centres hasn't been taken into account in this model.

4. Based on general knowledge and understanding of the authors with respect to land, electricity, water and other prices in Canada, the fixed costs and fixed carbon footprint associated with setting up of collection as well as reprocessing centres are randomly simulated within a given range of values.

5. Multiple reprocessing techniques are available for decontamination and sterilization of N95 masks such as use of vaporized hydrogen peroxide, ultraviolet germicidal radiation, ozone radiation, etc. [28]. However, when this study was carried out, there was a lack of sufficient literature mentioning the disinfection techniques for FFP2 respirators approved for use by Government of Canada since numerous developments, tests and trials were going on at that time. For running this model, we have taken ultraviolet germicidal radiation as the disinfection technique.

6. The selling price of a mask is taken as the price at which N95 masks were being sold on a renowned e-commerce platform for healthcare utilities, accessories and equipment in Canada [31]. Cost of disposal is assumed to be in the range of 20% to 30% of the price of the mask which accounts for the material costs.

7. Costs of collection and reprocessing are approximated in a particular range based on minimum hourly labour wages in Canada and costs of operating equipment related to ultraviolet germicidal radiation technique for decontamination of masks respectively.

8. Carbon footprint values associated with manufacturing and disposal of masks are taken for a N95 mask from literature [32].

9. Carbon footprint associated with reprocessing is taken to be equivalent to that of using ultraviolet germicidal radiation as the sterilization technique and that associated with collection process is assumed within a particular range of values based on understanding and experience of the authors.

## 4 Proposed Mathematical Model

The objective of this mathematical model for a closed-loop supply chain of PPE face masks is to optimize the reverse logistics network with respect to the location of collection, reprocessing, and recycling facilities. The goals are to maximize profits, minimize environmental impact in terms of carbon footprint, and maximize job creation.

### 4.1 Notations

We introduce the following notations:

**Indices:**

$i$ – Set of healthcare departments; 1, ..., $I$
$j$ – Set of collection centres; 1, ..., $J$
$k$ – Set of reprocessing centres; 1, ..., $K$
$l$ – Set of potential location sites; 1, ..., $L$
$m$ - Set of disposal sites; 1, ..., $M$

**Parameters:**

B – Maximum Budget allocated for setting up of collection centre $j$ and reprocessing centre $k$
$CC_j^l$ – Cost of collection at collection centre $j$ at location site $l$
$CD_m$ – Cost of disposal at disposal site $m$
$CO_2$ - Carbon Footprint emission of a small sized truck/km
CT – Cost of transportation in a small sized truck/km
$D_{i,j}^l$ – Distance between hospital $i$ and collection centre $j$
$D_{j,k}^l$ – Distance between collection centre $j$ and reprocessing centre $k$

E – Carbon footprint associated with production of a mask

$ECC_j^l$ – Carbon footprint during the collection process at collection centre $j$ at location site $l$

$ED_m$ – Carbon Footprint associated with disposal at disposal site $m$

$ERC_k^l$ – Carbon footprint during reprocessing at reprocessing centre k at reprocessing site $l$

$F_j^l$ – Fixed cost of opening collection centre $j$ at location $l$

$F_k^l$ – Fixed cost of opening reprocessing centre $k$ at location $l$

$FECC_j^l$ – Fixed carbon footprint associated with setting up of collection centre $j$ at location $l$

$FERC_k^l$ – Fixed carbon footprint associated with setting up of reprocessing centre k at location $l$

$JC_j^l$ – Number of jobs created at collection centre $j$ at location $l$

$JRC_k^l$ – Number of jobs created at reprocessing centre $k$ at location $l$

P – Selling price of a mask

$RC_k^l$ – Cost of reprocessing at reprocessing centre $k$ at location site $l$

$U_i^l$ – Usage of masks at hospital i

**Decision Variables:**

$Q_j$ – Number of masks at the collection centre $j$

$Q_k$ – Number of masks at the reprocessing centre $k$

$Q_m$ – Number of masks at the disposal centre $m$

$X_j^l \in \{0,1\}$ – Opening collection centre $j$ at location $l$

$X_k^l \in \{0,1\}$ – Opening reprocessing centre $k$ at location $l$

$Y_{i,j}^l \in \{0,1\}$ – If collection centre $j$ is assigned to demand from hospital $i$ at location $l$

$Y_{j,k}^l \in \{0,1\}$ – If reprocessing centre $k$ is assigned to collection centre $j$ at location $l$

## 4.2    Objective Functions

**Objective Function 1:**

The total costs can be broken down into fixed costs and variable costs. Fixed Costs include the cost associated with setting up of new facility like rent, basic utilities like water, electricity, etc. and are incurred irrespective of whether the facility if operational or not.

Fixed Costs $= \sum_j F_j^l X_j^l + \sum_k F_k^l X_k^l$ (1)

Variable Costs include transportation costs, operations costs at collection and reprocessing centres.

Variable Costs $= \sum_i \sum_j \mathrm{CT} * D_{i,j}{}^l Y_{i,j}{}^l + \sum_k \sum_j \mathrm{CT} * D_{j,k}{}^l Y_{j,k}{}^l + \sum_i \sum_j CC_j{}^l Q_j Y_{i,j}{}^l + \sum_k \sum_j RC_k{}^l Q_k Y_{j,k}{}^l + \sum_m \sum_j CD_m Q_m$ (2)

Total Revenue from the circular supply chain is the price of a mask at which it is sold multiplies by the total quantity of reprocessed masks.

Total Revenue $= \sum_k P * Q_k$ (3)

Since profit is equal to total revenue minus total costs, our first objective function is to maximize the profit (economic pillar of sustainability) in the supply chain.

Max $Z_1$ = Total Revenue – Fixed Costs – Variable Costs (4)

**Objective Function 2:**

Environmental Recovery due to reprocessing can be calculated as the carbon footprint associated with production of one mask multiplied by the total number of reprocessed masks.

Environmental Recovery $= \sum_k E * Q_k$ (5)

Carbon Footprint is generated due to transportation, operations at collection and reprocessing facilities and disposal of masks

Environmental Costs in terms of generated Carbon Footprint $= \sum_i \sum_j D_{i,j}{}^l CO_2 Y_{i,j}{}^l + \sum_k \sum_j D_{j,k}{}^l CO_2 Y_{j,k}{}^l + \sum_j FECC_j{}^l X_j{}^l + \sum_k FERC_k{}^l X_k{}^l + \sum_i \sum_j ECC_j{}^l Q_j Y_{i,j}{}^l + \sum_k \sum_j ERC_k{}^l Q_k Y_{j,k}{}^l + \sum_m \sum_j ED_m Q_m$ (6)

Our second objective function is to maximize net environmental advantage in this supply chain which is equal to environmental recovery minus environmental costs in terms of generated carbon footprint

Max $Z_2$ = Environmental Recovery – Environmental Costs in terms of generated carbon footprint (7)

**Objective Function 3:**

Maximize the number of jobs (social pillar of sustainability) created at the collection centre and the reprocessing centre.

Max $Z_3 = \sum_j JC_j{}^l X_j{}^l + \sum_k JRC_k{}^l X_k{}^l$ (8)

## 4.3    Constraints

Constraint 9 shows reprocessing centre k be selected if collection centre j associated with reprocessing centre k.

$X_k^l \geq \sum_j Y_{j,k}^l \quad \forall k \ (9)$

Constraint 10 shows collection centre j be selected if hospital i associated with collection centre j.

$X_j^l \geq \sum_i Y_{i,j}^l \quad \forall j \ (10)$

Constraint 11 shows that supply from hospitals needs to be supplied to at least one open collection centre j.

$\sum_j Y_{i,j}^l \geq 1 \quad \forall i \ (11)$

Constraint 12 shows that supply from collection centres needs to be supplied to at least one open reprocessing centre k.

$\sum_k Y_{j,k}^l \geq 1 \quad \forall j \ (12)$

Constraint 13 shows that the fixed costs associated with setting up of collection centre j and reprocessing centre k should be less than the total budget allocated.

$\sum_j F_j^l X_j^l + \sum_k F_k^l X_k^l \leq B \ (13)$

Constraint 14 shows that the total quantity of masks being collected at collection centre j is equal to alpha times the total quantity of masks at hospital i assigned to collection centre j.

$\sum_j \sum_i \alpha \times U_i \times Y_{i,j}^l = \sum_j Q_j \ (14)$

Constraint 15 shows that the total quantity of masks being collected at reprocessing centre k is equal to beta times the total quantity of masks at collection centre j assigned to reprocessing centre k.

$\sum_j \sum_k \beta \times Q_j \times Y_{j,k}^l = \sum_k Q_k \ (15)$

Constraint 16 shows that the total quantity of masks at reprocessing centre k subtracted from the total quantity of masks at collection centre j gives us the total number of masks sent for disposal to site m.

$\sum_j Q_j - \sum_k Q_k = \sum_m Q_m \ (16)$

## 5    Results and Discussions

The proposed circular economy model for PPE masks has been solved using the ε-constraint technique. This method, commonly used in multi-objective optimization problems, converts a multi-objective problem into a series of single-objective subproblems. In multi-objective optimization, conflicting objectives must be optimized simultaneously. Due to trade-offs between these objectives, finding a single solution that optimizes all objectives is often not feasible. The ε-constraint technique addresses this by breaking down the multi-objective problem into several single-objective subproblems. The mathematical model was coded and solved using the CPLEX solver in GAMS, running on an Intel Core i5 processor.

Table 1 compares the linear and proposed circular economy models for masks. The results highlighted significant improvements in the circular economy model across economic, environmental, and social dimensions. In particular, we came up the following observations:

- **Economic Perspective:** The circular economy model showed a profit of 38,992.187 CAD compared to the linear economy model, which resulted in a loss of 982,167.795 CAD. The loss in the linear model was calculated based on the assumption that all masks used at the hospital/healthcare center echelon are disposed of post-use, with the cost of disposal factored in.
- **Environmental Perspective:** The environmental impact was measured in terms of $CO_2$ emissions. The linear economy model generated approximately 65477.853 kg of $CO_2$, calculated by multiplying the total number of masks used by the carbon footprint associated with disposing each mask. In contrast, the circular economy model resulted in a net environmental benefit of 1766.334 kg of $CO_2$, demonstrating a significant reduction in environmental impact.
- **Social Perspective:** In terms of job creation, the circular economy model generated 53 jobs, compared to the 20-30 jobs created by the linear economy model.

**Table 1.** Comparison of proposed circular economy model and existing linear economy model for economic, environmental and social perspectives

| Model Type | $Z_1$ (CAD) | $Z_2$ (Kg of $CO_2$) | $Z_3$ (No. of jobs) |
|---|---|---|---|
| Circular Economy | 38992 | 1767 | 53 |
| Linear Economy | -982167 | -65477 | 20-30 |

10 more iterations were conducted, with the lower bound of Z2 and Z3 changed every time to simultaneously address all three objectives to determine ideal compromise solutions. The results of Pareto solutions given by applying ε-constraint technique and the corresponding Pareto chart is depicted in Figure 5.

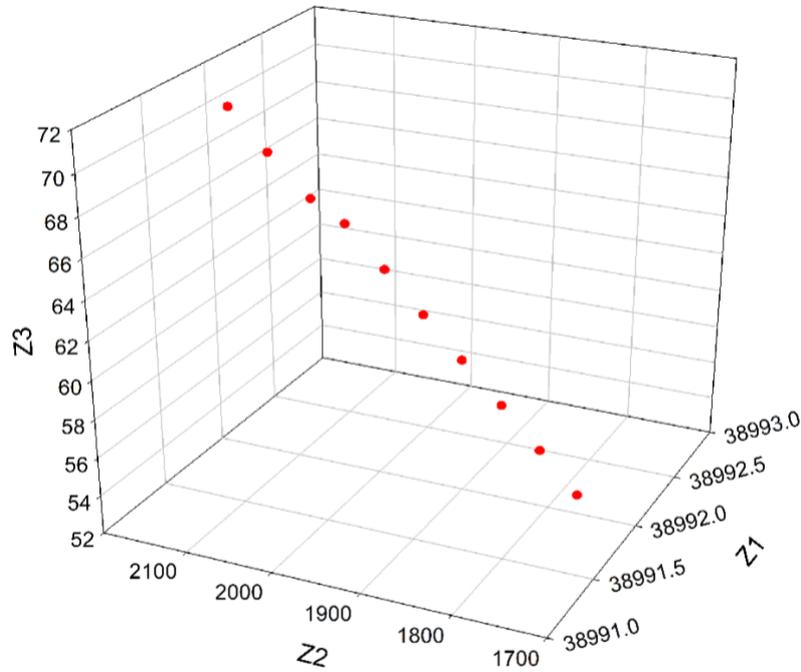

**Fig. 7.** Pareto front for economic, environmental, and social non-dominated solutions by ε-constraint method

From the potential location of facilities identified for setting up of collection and reprocessing facilities, twenty optimal locations for collection facilities as calculated by model were sites 1, 2, 3, 4, 5, 6, 8, 9, 10, 11, 12, 13, 14, 16, 17, 18, 20, 21, 23 and 24 whereas seven reprocessing facilities were 1, 5, 12, 14, 16, 18 and 24.

The optimization model calculated the optimal flow of masks from hospitals to collection centers and from collection centers to reprocessing centers or disposal sites. The results indicated that 7,275,317 masks are collected at the collection facilities from hospitals, of which 6,911,551 are sent to reprocessing facilities, and the remaining 363,766 are sent for disposal.

These results demonstrate the viability and advantages of the proposed circular economy model over the existing linear economy model. The circular model enhances economic profitability, significantly reduces environmental impact, and increases job creation, aligning with the three pillars of sustainability: economic, environmental, and social.

# 6    Conclusion and Future Work

This paper examined how Circular Economy (CE) model can improve the PPE supply chain, focusing on surgical masks and respirators in British Columbia, Canada. The goal was to tackle challenges of distribution uncertainty and demand surges during pandemics by creating a reverse logistics network. This network includes facilities for collecting, disinfecting, and reprocessing PPE with the aim of closing the loop on PPE use.

We developed a location-allocation optimization model using real-time data from September 2020. The model balanced economic, environmental, and social objectives, finding the optimal locations for collection and reprocessing centers from 25 potential sites. It also optimized the flow of masks through the reverse logistics network. To manage the conflicting objectives, we applied the ε-constraint method, which turned the multi-objective problem into a single objective by introducing constraints.

The results show that the CE model performs better than the LE model in terms of economic and social benefits, though environmental challenges remain. The CE model significantly reduced waste management costs and created local jobs, but the environmental gains were limited by higher energy and water usage. Despite these challenges, the model still reduced overall carbon emissions compared to the linear approach.

This model can be expanded to cover all of British Columbia and, eventually, the rest of Canada or any other country. It could also include more types of PPE, like gloves, shields, and gowns, turning it into a multi-product system. Future research could explore integrating reprocessing facilities within hospitals, where used PPE could be collected, disinfected, and reprocessed on-site. A cost-benefit analysis would help assess this option. Adding a recycling process to the CE framework could further cut down the use of raw materials, improving both environmental and economic sustainability.

The managerial implications include:

- The CE model reduces reliance on global suppliers, helping organizations manage supply chain disruptions more effectively.
- The CE model boosts profitability while also creating local jobs by bringing down waste management costs and using resources more efficiently.
- Moving to a CE model helps organizations lower their carbon emissions and create social value by integrating sustainability into their operations.
- The optimization model offers managers a clear strategy for selecting the best locations for collection and reprocessing centers, minimizing transportation costs and improving logistics.
- The CE model can be easily scaled from city to provincial or national levels and adapted for different contexts, thus making it highly flexible for various industries.
- A key takeaway is that supply chains must be designed with circularity in mind from the beginning and not as an afterthought. Building circular principles into the supply chain structure from the outset ensures smoother implementation, greater efficiency and long-term sustainability. This proactive approach is crucial to fully realizing the environmental, social, and economic benefits of a circular economy.

- Successfully implementing a circular economy model requires close collaboration with suppliers and other key stakeholders. Organizations need to work together to develop recyclable products and build the circular logistics infrastructure necessary to support collection, disinfection, reprocessing, and redistribution. Engaging stakeholders early in the process fosters innovation, aligns objectives, and ensures that the circular system is robust and effective.

The methodology and framework can be extended to other PPE components and industries, promoting a holistic approach to the adoption of the circular economy across different sectors. Engaging stakeholders, including suppliers, healthcare facilities and waste management entities is not just crucial but also integral for successfully implementing a circular economy model. Their collaborative efforts can drive innovation and shared value creation, making them an essential part of the process.

This approach can also be applied to other sectors to support a broader shift toward circular economy principles. By implementing the proposed circular economy model, managers can improve operational efficiency and profitability and contribute to broader sustainability goals. This will ensure a more resilient and sustainable supply chain for PPE and other critical healthcare supplies.

## References

1. Park, C., Kim, K., Roth, S., Beck, S., Kang, J. W., Tayag, M. C., & Grifin, M. (2020). *Global Shortage of Personal Protective Equipment amid COVID-19: Supply Chains, Bottlenecks, and Policy Implications*. https://doi.org/10.22617/brf200128-2
2. World Health Organization: WHO. (2020, March 3). *Shortage of personal protective equipment endangering health workers worldwide*. World Health Organization. https://www.who.int/news/item/03-03-2020-shortage-of-personal-protective-equipment-endangering-health-workers-worldwide
3. Yu, H., Sun, X., Solvang, W. D., & Zhao, X. (2020). Reverse Logistics Network Design for Effective Management of Medical Waste in Epidemic Outbreaks: Insights from the Coronavirus Disease 2019 (COVID-19) Outbreak in Wuhan (China). *International Journal of Environmental Research and Public Health*, *17*(5), 1770. https://doi.org/10.3390/ijerph17051770
4. MacNeill, A. J., Hopf, H., Khanuja, A., Alizamir, S., Bilec, M., Eckelman, M. J., Hernandez, L., McGain, F., Simonsen, K., Thiel, C., Young, S., Lagasse, R., & Sherman, J. D. (2020). Transforming the Medical Device Industry: Road map to a Circular economy. *Health Affairs*, *39*(12), 2088–2097. https://doi.org/10.1377/hlthaff.2020.01118
5. Kumar, H., Azad, A., Gupta, A., Sharma, J., Bherwani, H., Labhsetwar, N. K., & Kumar, R. (2020). COVID-19 Creating another problem? Sustainable solution for PPE disposal through LCA approach. *Environment Development and Sustainability*, *23*(6), 9418–9432. https://doi.org/10.1007/s10668-020-01033-0
6. Krysovatyy, A., Zvarych, I., & Zvarych, R. (2018). Circular economy in the context of alterglobalization. *JOURNAL OF INTERNATIONAL STUDIES*, *11*(4), 185–200. https://doi.org/10.14254/2071-8330.2018/11-4/13

7.  Vitacore Industries Inc. (2021, February 3). *Vitacore Launches Canada's First Recycling Program for Single-Use Masks and Respirators*. Cision News. https://www.newswire.ca/news-releases/vitacore-launches-canada-s-first-recycling-program-for-single-use-masks-and-respirators-891028844.html

8.  Wang, L., Li, S., Ahmad, I. M., Zhang, G., Sun, Y., Wang, Y., Sun, C., Jiang, C., Cui, P., & Li, D. (2023). Global face mask pollution: threats to the environment and wildlife, and potential solutions. *The Science of the Total Environment*, *887*, 164055. https://doi.org/10.1016/j.scitotenv.2023.164055

9.  Team, E. C. (2022, June 24). How can we recycle medical waste: COVID-19 masks› Evreka. *Evreka › How Can We Recycle Medical Waste: COVID-19 Masks*. https://evreka.co/blog/a-sustainable-way-to-handle-covid-19-masks/

10. Khan, M. T., Shah, I. A., Hossain, M. F., Akther, N., Zhou, Y., Khan, M. S., Al-Shaeli, M., Bacha, M. S., & Ihsanullah, I. (2022). Personal protective equipment (PPE) disposal during COVID-19: An emerging source of microplastic and microfiber pollution in the environment. *The Science of the Total Environment*, *860*, 160322. https://doi.org/10.1016/j.scitotenv.2022.160322

11. *Live Draw SDY: Live SDY hari ini, Togel Sidney 4D, Result SDY Pools tercepat*. (n.d.). Live Draw SDY. https://designedconscious.com/plastics-in-the-ocean/sustainability-news-stories/can-covid-19-masks-be-recycled/

12. Chen, C., Chen, J., Fang, R., Ye, F., Yang, Z., Wang, Z., Shi, F., & Tan, W. (2021). What medical waste management system may cope With COVID-19 pandemic: Lessons from Wuhan. *Resources Conservation and Recycling*, *170*, 105600. https://doi.org/10.1016/j.resconrec.2021.105600

13. Karliner, J., Slotterback, S., Boyd, R., Ashby, B., Steele, K., & Wang, J. (2020). Health care's climate footprint: the health sector contribution and opportunities for action. *European Journal of Public Health*, *30*(Supplement_5). https://doi.org/10.1093/eurpub/ckaa165.843

14. Betcheva, L., Erhun, F., & Jiang, H. (2020). OM Forum—Supply Chain Thinking in Healthcare: Lessons and Outlooks. *Manufacturing & Service Operations Management*, *23*(6), 1333–1353. https://doi.org/10.1287/msom.2020.0920

15. Wieser, P. (2011). From health logistics to health supply chain management. *Supply Chain Forum an International Journal*, *12*(1), 4–13. https://doi.org/10.1080/16258312.2011.11517249

16. Gendy, A. W. A., & Lahmar, A. (2019). Review on Healthcare Supply Chain. In *ACS/IEEE International Conference on Computer Systems and Applications* (pp. 1–10). https://doi.org/10.1109/aiccsa47632.2019.9035234

17. Mittal, A., & Mantri, A. (2023). A literature survey on healthcare supply chain management. *F1000Research*, *12*, 759. https://doi.org/10.12688/f1000research.131440.2

18. MacNeill, A. J., Hopf, H., Khanuja, A., Alizamir, S., Bilec, M., Eckelman, M. J., Hernandez, L., McGain, F., Simonsen, K., Thiel, C., Young, S., Lagasse, R., & Sherman, J. D. (2020b). Transforming the Medical Device Industry: Road map to a Circular economy. *Health Affairs*, *39*(12), 2088–2097. https://doi.org/10.1377/hlthaff.2020.01118

19. Geissdoerfer, M., Savaget, P., Bocken, N. M., & Hultink, E. J. (2016). The Circular Economy – A new sustainability paradigm? *Journal of Cleaner Production*, *143*, 757–768. https://doi.org/10.1016/j.jclepro.2016.12.048

20. Chakravarty, A. K. (2013). Sustainable supply chains. In *Springer texts in business and economics* (pp. 273–305). https://doi.org/10.1007/978-3-642-41911-9_9

21. Sambrani, V., & Pol, D. (2016, October 16). *An overview of sustainable supply chain management*. https://papers.ssrn.com/sol3/papers.cfm?abstract_id=2870430

22.  Bouchery, Y., Corbett, C. J., Fransoo, J. C., & Tan, T. (2016). Sustainable supply chains. In *Springer series in supply chain management*. https://doi.org/10.1007/978-3-319-29791-0

23.  Formentini, M. (2020). Sustainable Supply chain management. In *Management for professionals* (pp. 207–223). https://doi.org/10.1007/978-3-030-56344-8_12

24.  Medhekar, A. (2023). Sustainable Supply Chains for Circular Economy in the Health sector: Challenges and Opportunities Post Pandemic. In *Advances in finance, accounting, and economics book series* (pp. 429–448). https://doi.org/10.4018/978-1-6684-7664-2.ch021

25.  Vishwakarma, A., Dangayach, G., Meena, M., Gupta, S., Joshi, D., & Jagtap, S. (2022). Can circular healthcare economy be achieved through implementation of sustainable healthcare supply chain practices? Empirical evidence from Indian healthcare sector. *Journal of Global Operations and Strategic Sourcing*, *17*(2), 230–246. https://doi.org/10.1108/jgoss-07-2022-0084

26.  Senna, P., Reis, A., Marujo, L. G., De Guimarães, J. C. F., Severo, E. A., & De Souza Gomes Dos Santos, A. C. (2023). The influence of supply chain risk management in healthcare supply chains performance. *Production Planning & Control*, *35*(12), 1368–1383. https://doi.org/10.1080/09537287.2023.2182726

27.  *Towards the circular economy Vol. 1: an economic and business rationale for an accelerated transition*. (2013, January 1). https://www.ellenmacarthurfoundation.org/towards-the-circular-economy-vol-1-an-economic-and-business-rationale-for-an

28.  Canadian Coalition for Green Health Care. (2024, June 18). *PPE-MSUP Project - Canadian Coalition for Green Health Care Inc.* Canadian Coalition for Green Health Care Inc. - the Canadian Coalition for Green Health Care. https://greenhealthcare.ca/ppe-msup/

29.  Yaspal, B., Jauhar, S. K., Kamble, S., Belhadi, A., & Tiwari, S. (2023). A data-driven digital transformation approach for reverse logistics optimization in a medical waste management system. *Journal of Cleaner Production*, *430*, 139703. https://doi.org/10.1016/j.jclepro.2023.139703

30.  Yadav, S., Mangla, S. K., Priyamvada, P., Borkar, A., & Khanna, A. (2023). An energy-efficient model for PPE waste management in a closed-loop supply chain. *Business Strategy and the Environment*, *33*(2), 1191–1207. https://doi.org/10.1002/bse.3541

31.  *N95 Mask, N95 Mask Canada, and N95 Respirator*. (n.d.). 72hours.ca. https://72hours.ca/collections/n95-mask

32.  Luo, Y., Yu, M., Wu, X., Ding, X., & Wang, L. (2023). Carbon footprint assessment of face masks in the context of the COVID-19 pandemic: Based on different protective performance and applicable scenarios. *Journal of Cleaner Production*, *387*, 135854. https://doi.org/10.1016/j.jclepro.2023.135854

# Acknowledgements

This research was carried out by Jainil Dharmil Shah under the guidance and supervision of Dr. Adel Guitouni (Associate Professor), Dr. Behrooz Khorshidvand (Post Doctoral Research Fellow) & Niloofar Gilani Larimi (PhD Candidate) at the Peter B. Gustavson School of Business, University of Victoria, Victoria, BC, Canada. Jainil was then a visiting undergraduate summer researcher at the University of Victoria as a part of Mitacs Globalink Research Internship Program '23 and is currently enrolled for his graduate studies at the Edwardson School of Industrial Engineering, Purdue University, West Lafayette, IN, USA starting from Fall '24.